\documentclass{amsart}

\usepackage{amsmath}
\usepackage{amssymb}
\usepackage{mathrsfs}
\usepackage{url}
\usepackage{enumitem}
\usepackage{stmaryrd}
\usepackage[all]{xy}
\usepackage{verbatim}

\newtheorem{Thm}{Theorem}[section]
\newtheorem{Con}[Thm]{Conjecture}

\newtheorem{Def}[Thm]{Definition}

\numberwithin{equation}{Thm}

\newcommand{\CC}{\mathbb{C}}

\newcommand{\FF}{\mathbb{F}}

\newcommand{\QQ}{\mathbb{Q}}

\newcommand{\ZZ}{\mathbb{Z}}

\date{\today}

\author{Daqing Wan}
\address{Department of
Mathematics, University of California, Irvine, CA 92697}
\email{dwan@math.uci.edu}

\title{Class numbers and $p$-ranks in $\ZZ_p^d$-towers}

\begin{document}

\begin{abstract}
To extend Iwasawa's classical theorem from $\ZZ_p$-towers to $\ZZ_p^d$-towers, Greenberg conjectured that 
the exponent of $p$ in the $n$-th class number in a $\ZZ_p^d$-tower of a global field $K$ ramified at finitely many primes is given by a polynomial 
in $p^n$ and $n$ of total degree at most $d$ for sufficiently large $n$. This conjecture remains open for $d\geq 2$. In this paper, we prove that this conjecture 
is true in the function field case.  Further, we propose a series of general conjectures on $p$-adic stability of zeta functions in a $p$-adic Lie tower of function fields. 
\end{abstract}

\maketitle

%\section{Introduction}

\section{Introduction} 

Let $K$ be a global field, that is, either a number field or a global function 
field with constant field $\FF_q$, where $\FF_q$ is the finite field of $q$ elements with characteristic $p$. 
In the second case, $K$ is simply the function field of a smooth projective geometrically 
irreducible curve $C$ defined over $\FF_q$. 
For a positive integer $d$,  
consider a $\ZZ_p^d$-tower of global fields 
$$K_{\infty}=\bigcup_{n=0}^{\infty} K_n \supset \cdots  \supset K_n \supset \cdots \supset K_1 \supset K_0=K,$$
where for all $n\geq 0$, 
$${\rm Gal}(K_n/K)  = (\ZZ_p/p^n\ZZ_p)^d.$$
In this paper, we always assume that the tower is ramified only at finitely many primes of $K$. 
This condition is automatically satisfied in the number field case.  The condition is natural and necessary in the function field case 
in order for many definitions and theorems to make sense. Let $h_n$ denote the class number of the field $K_n$.  
The basic problem in Iwasawa theory is to understand the stable behavior of the $p$-part of $h_n$, namely, the $p$-adic 
valuation $v_p(h_n)$ as a function of $n$.  In the literature, the number $v_p(h_n)$ is usually denoted by $e_n$, which is 
the exponent of $p$ in the class number $h_n$. 
To extend Iwasawa's classical result for $\ZZ_p$-extensions, 
Ralph Greenberg made the following conjecture about $40$ years ago, see section $7$ in 
\cite{CM}. 
 
\begin{Con}\label{Con1.0} Let $K_{\infty}/K$ be a $\ZZ_p^d$-extension of a global field $K$ which ramifies at finitely 
many primes. 
There is a polynomial $E(x,y)\in  \QQ[x,y]$ of total degree at most $d$ depending on the tower such that for all sufficiently large $n$, we have 
$$v_p(h_n) = E(p^n, n).$$
\end{Con}

In the case $d=1$ and $K$ is a number field, this is precisely the classical Iwasawa theorem (1959) which says that 
$v_p(h_n) = \mu p^n +\lambda n +\nu$ for all sufficiently large $n$, where $\mu, \lambda, \nu$ are constants depending on the tower. 
Historically, Iwasawa first observed that in the case $d=1$ and $K_{\infty}$ is a constant field $\ZZ_p$-extension of a function field, 
the statement is true. This motivated him to develop the theory of $\ZZ_p$ extensions, which led to his proof of the case $d=1$ for number fields. 
It was realized in Mazur-Wiles \cite{MW} (1983) and Gold-Kisilevsky \cite{GK} (1988) that Iwasawa's ideas carry over to the function field case as well for $d=1$, not necessarily constant field 
extensions. Thus, Greenberg's conjecture 
was long known to be true in the case $d=1$, for both number fields and function fields. 

To understand how the $\mu$ and $\lambda$ invariants vary as $K_{\infty}$ varies over all $\ZZ_p$-extensions of a fixed number field $K$, 
Greenberg \cite{Gr} (1973) initiated the study of Iwasawa theory for $\ZZ_p^d$-extensions of  $K$. 
The subsequent development of $\ZZ_p^d$-extensions led Monsky \cite{Mo2} to prove that the $\mu$ invariant $\mu(K_{\infty}/K)$ is 
bounded by a constant depending only on $K$ and $p$, as $K_{\infty}$  varies over all $\ZZ_p$-extensions of a number field $K$. 
It seems still open if the $\lambda$ invariant $\lambda(K_{\infty}/K)$ is also bounded by a constant depending only 
on $K$ and $p$. 
Recall that for a number field $K$, the 
composite of all $\ZZ_p$-extensions of $K$ is a $\ZZ_p^r$-extension of $K$, where $r_2+1 \leq r \leq [K:\QQ]$ and $r_2$ is the number of complex 
primes of $K$. Leopoldt's conjecture is equivalent to saying that $r= r_2 +1$. In contrast, 
for a function field $K$, the composite of all $\ZZ_p$-extensions of $K$ is a huge $\ZZ_p^{\infty}$-extension, and the $\mu$-invariant 
$\mu(K_{\infty}/K)$ is not known to be bounded as $K_{\infty}$  varies over all $\ZZ_p$-extensions of $K$, see [LZ]. 

%It is not hard to prove  that the $\lambda$-invariant $\lambda(K_{\infty}/K) $ is unbounded for a fixed base function field $K$ by increasing the ramification locus. 

%Furthermore, both the $\mu$ and $\lambda$-invariants are bounded when 
%$K_{\infty}$ varies over all $\ZZ_p$-subextensions in any fixed 
%$\ZZ_p^d$-tower over $K$.  

In the number field case with $d\geq 2$, Greenberg's conjecture was initially studied by Cuoco \cite{Cu} (1980) in the case $d=2$, 
and later in a series of papers 
by Cuoco-Monsky \cite{CM} (1981), and Monsky \cite{Mo1}\cite{Mo2} \cite{Mo3}\cite{Mo4} (1981-1989) for general $d$. The strongest result so far is from Monsky \cite{Mo4} (1989) who proves the number field case of the following asymptotic theorem.  

\begin{Thm}\label{Monsky} Let $K_{\infty}/K$ be a $\ZZ_p^d$-extension of a global field $K$ which ramifies at finitely 
many primes. There are integers $m_0$, $\ell_0$ and a real number $\alpha$ depending on the tower such that 
for all sufficiently large $n$, we have $$v_p(h_n) = (m_0p^n +\ell_0n +\alpha)p^{(d-1)n} +O(np^{(d-2)n}).$$
\end{Thm}

As later remarked in Li-Zhao \cite{LZ} (1997), 
Monsky's proof carries over to function fields too, and thus the above asymptotic theorem is true for both number fields and function fields. 
Since then, progress on Greenberg's conjecture has stopped in both the number field case and the function field case, 
perhaps due to the following two reasons.  First, Cuoco-Monsky \cite{CM} wrote ``Despite the
evidence of this paper and of [4] we believe this (Greenberg conjecture) to be false, and show that a related
module-theoretic conjecture fails, even when $d = 2$".  In another paper, Monsky \cite{Mo1} wrote ``Greenberg's conjecture is 
probably false; the last section of [2] presents module-theoretic evidence against it".  
Second, Li-Zhao wrote 
``Classically, there are several approaches to this problem: (A) via class
number formula and $p$-adic L-functions; (B) via class number formula and
$p$-adic measures on $\ZZ_p$; (C) via the theory of noetherian modules over Iwasawa
algebra. Approach (A) works particularly well for the so called cyclotomic $\ZZ_p$-extensions. 
But it apparently fails for more general fields which do not possess a good class number formula. (B) was found by
Sinnott  when he studied the work of Ferrero and Washington.
Neither (A) nor (B) can be used to handle the function field case
because good analogues of $p$-adic L-functions and $p$-adic measures are still unborn for function fields. 
Approach (C) was the one initiated by Iwasawa and further developed by Serre, Greenberg, Cuoco, Monsky and many others.
This method can be combined with class field theory and Kummer theory to
get a lot of nice results. The last approach of the above also has the merit of being generalizable
to function field case, which is what we shall use in the present paper". This was the status up to 1997, and no further work 
was done on this conjecture during the last twenty years. 

Our aim of this paper is to prove that Greenberg's conjecture is true in the function field case for all $d\geq 1$. Namely, we have 

\begin{Thm}\label{Wan} Let $K_{\infty}/K$ be a $\ZZ_p^d$-tower of function fields of characteristic $p$ which ramifies at finitely 
many primes. Then, there is a polynomial $E(x,y)\in  \QQ[x,y]$ of total degree at most $d$ (and degree at most $1$ in $y$) 
depending on the tower such that for all sufficiently large $n$, 
we have 
$$v_p(h_n) = E(p^n, n).$$
\end{Thm}

We use the class number formula and $p$-adic L-function approach. 
The required universal $p$-adic L-function in function fields is the $T$-adic L-function introduced in [Liu-Wan] (2009) for $d=1$ and 
again in [RWXY] (2017) for $d>1$. The $T$-adic L-function is instrumental in the study of slope stability of zeta functions in 
$\ZZ_p^d$-tower, which is a new direction, see \cite{DWX, RWXY}. 
Its analytic property on the closed unit disk (the main conjecture) follows from the work of Crew \cite{crew} (1987) 
%and Emerton-Kisin \cite{EK} (2001) 
on a conjecture of Katz \cite{Ka}. The class number is a special value of the zeta function, which can be decomposed 
in terms of L-functions of the finite characters of the Galois group. 
To complete the proof, we partition the finite character space of the Galois group into a disjoint union of subspaces in terms of the exact ramification locus 
of the character. These disjoint subspaces turn out to be semi-algebraic in the sense of \cite{Mo1}, and one can then apply  Monsky's crucial power series 
lemma \cite{Mo1}. 
Contrary to Monsky's belief, philosophically our result suggests that Greenberg's conjecture should be true 
in the number field case. 

A bye-product of our proof also gives a stable formula for the geometric $p$-rank $r_p(n)$ of $K_n$. 
Recall that $K_n$ is the function field of a smooth projective geometrically irreducible curve $C_n$ defined 
over some finite extension field over $\FF_q$.  The Jacobian variety $J_n$ of $C_n$ is an abelian 
variety of dimension $g_n$, where $g_n$ is the genus of $C_n$. The $p$-adic Tate module is 
$$T_p(J_n) = \lim_{\leftarrow~k} J_n [p^k],$$
where $J_n[p^k]$ denotes the group of $p^k$-torsion points of $J_n$ over an algebraic closure of $\FF_q$. 
The Tate module $T_p(J_n)$ is a free $\ZZ_p$-module whose rank is 
called the $p$-rank of $K_n$, denoted by $r_p(n)$. 
It is well known that $0\leq r_p(n)\leq g_n$, where $g_n$ is the genus of $K_n$.  

\begin{Thm}\label{rank} Let $K_{\infty}/K$ be a $\ZZ_p^d$-tower of function fields of characteristic $p$ which ramifies at finitely 
many primes. Then, there is a polynomial $R(x) \in  \ZZ[x]$ of degree at most $d$
depending on the tower such that for all sufficiently large $n$, 
we have 
$$r_p(n) = R(p^n).$$
\end{Thm}

Finally, we note that there is no stable formula for the genus sequence $g_n$ in the function field case, even for $d=1$. 
In fact, the genus $g_n$ can grow as fast as one 
wants as a function of $n$. A simple necessary and sufficient condition for the stability of $g_n$  is given in \cite{KW} in the additive 
framework of Artin-Schreier-Witt construction of $\ZZ_p$-towers. For natural $\ZZ_p^d$-towers (those coming from algebraic geometry), we conjecture 
that the genus $g_n$ is  periodically given by  several polynomials in $p^n$ for all sufficiently large $n$. This is a small part of a 
series of general conjectures on $p$-adic stability of zeta functions for $p$-adic Lie towers of function fields that we shall state and discuss 
in the last section. It is hoped that these conjectures would motivate further work  in both the function field case and the number field case.

\section{The case with no constant $\ZZ_p$-extension}

Let $K$ be a global function field with constant field $\FF_q$. 
In this section and next section, we assume that the $\ZZ_p^d$-tower $K_{\infty}/K$ contains no constant subextension. 
The tower gives a continuous group isomorphism 
$$\rho: G_{\infty}:= {\text{Gal}}(K_{\infty}/K) \cong \ZZ_p^d.$$
For each integer $d$-tuple $(n_1,..., n_d) \in \ZZ_{\ge 0}^d$,  reduction modulo $(p^{n_1}, ..., p^{n_d})$ produces a 
subextension $K_{n_1,..., n_d}$ such that 
$$ G_{n_1,..., n_d}:= {\text{Gal}}(K_{n_1,..., n_d}/K) \cong \prod_{i=1}^d \ZZ_p/{p^{n_i}\ZZ_p}.$$
This makes sense even if some of the $n_i$ are $\infty$. 
For $1\leq i\leq d$, let $K_n^{(i)} = K_{0,.., n,..., 0}$, where $n$ lies in the $i$-th coordinate. For each $i$,  
$K_{\infty}^{(i)}/K$ is a $\ZZ_p$-tower and thus it is totally ramified at a prime $x$ if the first layer $K_1^{(i)}$ is ramified at $x$.  
It is clear that for $0\leq n \leq \infty$, 
$$K_n = K_{n,..., n} = K_{n,0,..., 0}\cdots K_{0,..., 0, n} =K_n^{(1)} \cdots K_n^{(d)}.$$
Let $P$ be the set of primes of $K$. 
Let $U$ be the unramified locus of the tower, which is a subset of $P$. 
The tower is ramified on $P-U$. By our assumption, $P-U$ is finite. Furthermore, by class field theory, 
the set $P-U$ is non-empty since there is no constant subextension.  
Each $p$-adic character $\chi: G_{\infty} \rightarrow \CC_p^*$ has a unique decomposition  
$$\chi = \chi_1\otimes \cdots \otimes \chi_d, $$
where each $\chi_i$ is a $p$-adic character of $G_{\infty}$ factoring through the quotient 
$$G_{\infty}^{(i)} = {\rm{Gal}}(K_{\infty}^{(i)}/K) \cong \ZZ_p.$$ 
To recall the class number formula, we need to introduce the zeta function.

Let $P_n$ denote the set of primes of $K_n$. 
Recall that the zeta function of $K_n$ is defined by 
$$Z(K_n,s) = \prod_{x\in P_n} \frac{1}{1-s^{{\rm deg}(x)}} \in 1 +s\ZZ [[s]].$$
The Riemann-Roch theorem implies that the zeta function is a rational function in $s$ of the form 
$$Z(K_n, s) = \frac{P(K_n, s)}{(1-s)(1-qs)}, \ P(K_n, s) \in 1 +s\ZZ[s],$$
where the numerator $P(K_n, s)$ is a polynomial in $s$ of degree $2g_n$ and $g_n= g(K_n)$ denotes the 
genus of $K_n$. 
By the celebrated theorem of Weil, the polynomial $P(K_n, s)$ is 
pure of $q$-weight $1$, that is, the reciprocal roots of $P(K_n, s)$ all have complex absolute value equal to $\sqrt{q}$. 
The class number $h_n$ of $K_n$ is given by the residue formula 
$$h_n= P(K_n, 1).$$
The Weil bound implies a good estimate for the complex absolute value of the class number $h_n$: 
$$(\sqrt{q}-1)^{2g_n} \leq h_n \leq (\sqrt{q}+1)^{2g_n}.$$
Thus, the growth of $h_n$ depends very much on the growth of the genus $g_n$ which can be complicated 
in general, unless the genus sequence $g_n$ becomes stable in some sense, see the last section for the genus stability conjecture for towers 
coming from algebraic geometry.  
In this paper, we shall mainly be interested in the $p$-adic absolute value of the class number, equivalently, $v_p(h_n)$ as a function of $n$.  
For this purpose, we use L-functions to decompose it. 

Recall that for positive integer $n\geq 1$, the Galois group 
$$G_{n}= {\text{Gal}}(K_{n}/K) \cong (\ZZ/{p^n\ZZ})^d.$$
For a finite continuous $p$-adic character $\chi: G_{\infty} \rightarrow \CC_p^*$, $\chi$ will factor through $G_n$ for some finite $n$. 
That is,  $\chi: G_{\infty} \rightarrow G_n  \rightarrow \CC_p^*$. The  L-function of $\chi$ over $K$ is 
$$L(\chi, s) = \prod_{x\in P, \chi ~{\rm{unramified ~ at}}~ x} \frac{1}{1-\chi(\text{Frob}_x)s^{\text{deg}(x)}} \in 1 +s\CC_p[[s]],$$
where  $\text{Frob}_x$ denotes the arithmetic Frobenius 
element of $G_{\infty}$   at $x$. Note that the ramification locus of the character $\chi$ could be smaller than $P-U$.  
If $\chi=1$, it is unramified everywhere and the L-function $L(\chi, s)$ is just the zeta function of $K$. If $\chi\not=1$, 
the L-function $L(\chi, s)$ is a polynomial 
in $s$, pure of weight $1$. 
One has the decomposition 
$$Z(K_n, s) = \prod_{\chi: G_{\infty}\rightarrow G_n \rightarrow \CC_p^*} L(\chi, s), \ 
\frac{P(K_n, s)}{P(K_0,s)} = \prod_{\chi\not= 1: G_{\infty} \rightarrow G_n \rightarrow \CC_p^*} L(\chi, s).$$
In order to understand each L-function $L(\chi, s)$, we need to know the ramification 
information of each finite character $\chi$. 

Let $W$ denote the group of all $p$-power roots of unity in $\mathbb{C}_p$. 
Let $X$ denote the group of finite continuous $p$-adic characters $\chi: G_{\infty} \rightarrow W$. 
There is a perfect pairing of $\mathbb{Z}_p$-modules 
$$\mathbb{Z}_p^d \times  W^d \longrightarrow W,\ ((a_1,\cdots, a_d), (\eta_1,\cdots, \eta_d)) \longrightarrow \prod_{i=1}^d \eta_i^{a_i}.$$
It follows that the group $X$ is isomorphic to $W^d$. Monsky \cite{Mo1} defines 
a Noetherian topology on $X=W^d$ with the property that the closed subsets of $W^d$ 
are finite unions of the basic set of the form 
$$\{ (\eta_1,\cdots, \eta_d) \in W^d |\prod_{i=1}^d \eta_i^{a_i} =\eta\},$$ 
where $\eta \in W$ and $(a_1,\cdots, a_d) \in \mathbb{Z}_p^d$. The open sets 
in $W^d$ are simply the compliments of closed subsets in $W^d$. Open sets and 
closed sets in $W^d$ are special cases of a more general class of sets, 
called semi-algebraic subsets in $W^d$, see \cite{Mo1}.

Next, we need to partition the character group $X$ in terms of their exact ramification locus. 
For each $x \in P-U$, let $I_x$ denote the inertial subgroup of $G_{\infty}$, which is a non-trivial $p$-adic subgroup of $G_{\infty}$. 
Let $X^*$ denote the ``interior piece" of $X$ consisting of those characters $\chi$ in $X$ 
such that $\chi (I_x)\not= 1$ for all $x \in P-U$. For such an interior character $\chi \in X^*$, it is clear that the ramification locus for $\chi$ is exactly 
$P-U$. Since $I_x$ is a finitely generated $\mathbb{Z}_p$-submodule of 
$\mathbb{Z}_p^d$, it follows that $X^*$ is an open subset of $X$.  

More generally, for any subset $S\subseteq P-U$, let $X_S$ denote the set of characters $\chi \in X$ such that 
$\chi(I_x)=1$ for all $x\not\in S$. Let $X_S^*$ denote the interior piece of $X_S$ consisting of those characters 
$\chi \in X_S$ such that $\chi(I_x) \not=1$ for all $x\in S$. Then, the ramification locus of $\chi \in X_S^*$ is 
exactly $S$, and we obtain the 
following disjoint interior decomposition 
$$X= \bigsqcup_{S \subseteq P-U} X_{S}^*, \ X^* = X_{P-U}^*.$$

For a positive integer $n$, let $X_n$ (resp., $X_{n, S}^*$) denote the set of characters $\chi \in X$ (resp., $\chi \in X_S^*$) such that 
$\chi^{p^n} =1$. Similarly, we have 
the 
following disjoint interior decomposition 
$$X_n= \bigsqcup_{S \subseteq P-U} X_{n, S}^*, \ X_n^* = X_{n, P-U}^*.$$
The class number $h_n$ then has the interior decomposition
$$\frac{h_n}{h_0} = \prod_{\chi\in X_n, \chi\not=1} L(\chi, 1) = \prod_{S\not= \phi, S \subseteq P-U} \prod_{\chi \in X_{n,S}^*} L(\chi, 1).$$

For any subset $S\subseteq P-U$, let $I_S$ denote the closure of subgroups $I_x$ for all $x\not\in S$. By the structure theorem 
for finitely generated $\mathbb{Z}_p$-modules, the quotient $G_{\infty}/I_S$ is of the form 
$$G_{\infty}/I_S \cong \mathbb{Z}_p^{d(S)} \oplus {H_S},$$
where $H_S$ is a finite abelian $p$-group and $0\leq d(S) \leq d$. Replacing the base field $K$ by $K_m$ for some large $m$ if necessary, we 
may assume that $H_S=0$ for all $S\subseteq P-U$. Thus, we have 
$$G_{\infty}/I_S \cong \mathbb{Z}_p^{d(S)}.$$
The characters in $X_S^*$ are exactly the interior piece in the $\mathbb{Z}_p^{d(S)}$-tower $K_{\infty, S}$ of $K$ arising from the fixed field in $K_{\infty}$ of $I_{P-U-S}$. 
Thus, without loss of generality, we can assume that $S=P-U$ and $X_{n,S}^* = X_n^*$. 

The above interior decomposition reduces Greenberg's conjecture for the polynomial stability of $v_p(h_n)$  to 
the following theorem. 

\begin{Thm}\label{Mon} 
There is a polynomial $E(x,y)\in  \QQ[x,y]$ of total degree at most $d$ (and degree at most $1$ in $y$)  such that for all sufficiently large $n$, 
we have 
$$v_p(\prod_{\chi \in X_{n}^*} L(\chi, 1)) = E(p^n, n).$$
\end{Thm}

Recall that the $p$-rank $r_p(n)$ of $K_n$ is well known to be the number of $p$-adic unit roots of $P(K_n, s)$. 
For $\chi\not = 1$, let $\ell_p(\chi)$ denote the number of $p$-adic unit roots for $L(\chi, s)$. The interior decomposition of the 
character space $X_n$ gives the following interior decomposition of the $p$-rank: 
$$r_p(n) - r_p(0) = \sum_{S\not= \phi, S \subseteq P-U} \sum_{\chi \in X_{n,S}^*} \ell_p(\chi).$$

The set $X_S^*$ is an open subset of $X_S = W^{d(S)}$. 
Lemma 4.1 (or Theorem 4.6) in Monsky \cite{Mo1} implies that the cardinality of 
the set $X_{n,S}^*$ is a polynomial in $p^n$ of degree bounded by $d(S)$ for all large $n$. To prove Theorem \ref{rank}, it is enough to prove the 
following theorem. 

\begin{Thm}\label{rank1} 
There is a non-negative integer $c$ depending on the tower such that for all sufficiently large $n$ and all $\chi \in X_n^*$, 
we have 
$$\ell_p(\chi) = c.$$
\end{Thm}

To prove the above two theorems, we need to package all the L-functions $L(\chi, s)$ for all characters $\chi \in X^*$ into 
a single universal $p$-adic L-function and study its variation as $\chi$ varies. This is accomplished by the $T$-adic L-function that we shall introduce next.

\section{$T$-adic L-functions} 

We now define the $T$-adic L-function first introduced in [LW] for $d=1$ and later in [RWXY] for general $d$.  
Instead of just finite characters $\chi_n: G_{\infty} \rightarrow G_n \rightarrow \CC_p^*$, we will also consider all continuous $p$-adic characters $\chi: G_{\infty} \longrightarrow \CC_p^*$, not necessarily of finite order. 
The isomorphism 
$$\rho: G_{\infty} \cong \ZZ_p^d$$ 
is crucial for us.  Let $\ZZ_p[[T]]=\ZZ_p[[T_1,\cdots, T_d]]$ be the power series ring in $d$ variables over $\ZZ_p$. 
Consider 
the universal continuous $T$-adic character 
$$\ZZ_p^d \longrightarrow \ZZ_p[[T]]^*, \  (a_1, \cdots, a_d) \longrightarrow \prod_{i=1}^d(1+T_i)^{a_i} \in \ZZ_p[[T]]^*. $$
Composing this universal $T$-adic character of $\ZZ_p^d$ with the isomorphism $\rho$, we get the universal $T$-adic character of $G_{\infty}$: 
\[
\rho_T: G_{\infty} \longrightarrow \ZZ_p^d \longrightarrow \text{GL}_1(\ZZ_p[[T]]) = \ZZ_p[[T]]^*.
\]
Let $D_p(1)$ denote the open unit disk $\{t=(t_1,\cdots, t_d)| |t_i|_p<1\}$ in $\CC_p^d$. 
For any element $t=(t_1,\cdots, t_d)  \in D_p(1)$, we have a natural evaluation map $\ZZ_p[[T]]^* \rightarrow \CC_p^*$ sending $T=(T_1, \cdots, T_d) $ 
to $t=(t_1, \cdots, t_d)$.  Composing all these maps, we get, for fixed $t\in D_p(1)$, a 
continuous character
\begin{equation} \label{varying character}
\rho_t: G_{\infty} \longrightarrow \CC_p^*. 
\end{equation}
The open unit disk $D_p(1)$ parametrizes all continuous $\CC_p$-valued  characters $\chi=\chi_1\otimes \cdots \otimes \chi_d$ of $G_{\infty}$ via the relation 
$t = t_{\chi}= (\chi_1(1)-1, \cdots, \chi_d(1)-1)$, where $\chi_i$ factors through the quotient $G_{\infty}^{(i)}: = {\rm Gal}(K_{\infty}^{(i)}/K) \cong \ZZ_p$. 
The L-function of $\rho_t$ is defined in the usual way: 
$$L({\rho_t}, s) = \prod_{x\in P, \rho_t ~{\rm{unramified ~ at}}~ x}\frac{1}{1-\rho_t({\rm Frob}_x)s^{{\rm deg}(x)}}\in 1 +s\CC_p[[s]].$$
For a general $t \in D_p(1)$, this L-function $L(\rho_t, s)$ does not  have a $p$-adic meromorphic continuation to the $p$-adic plane 
$|s|_p <\infty$, see \cite{Wa0} for results in this direction. 

In the case that $\chi=\chi_1\otimes \cdots \otimes \chi_d$ is a finite $p$-adic character,  
we have 
$$t=(\chi_1(1)-1,\cdots, \chi_d(1)-1) \ , \ L(\chi, s) = L({\rho}_{t}, s). $$ 
Elements of the form $t_{\chi}=(\chi_1(1)-1, \cdots, \chi_d(1)-1)$ with each $\chi_i$ finite order,  are called the classical points in $D_p(1)$. 
The classical points are exactly the elements in 
$$W^d - 1 = \{(\eta_1-1, \cdots, \eta_d-1) | \eta_i \in W \}.$$
For a classical point $t \in W^d-1$, the ramification locus of the corresponding character $\rho_t \in X$ may be strictly smaller than $P-U$. 
Let $D^*$ be the set of classical points $t\in W^d-1$ such that 
the character $\rho_t \in X^*$ (the interior characters). We call $D^*$ as the 
set of interior classical points. 
Then, the ramification locus of $\rho_t$ 
is exactly $P-U$ for all interior classical points $t \in D^*$. 

As the order of $\chi$ goes to infinity, $t_{\chi}$ approaches to the boundary of the disk $D_p(1)$. 
Thus, to understand the behavior of $L(\chi, s)$ as ${\rm ord}(\chi)$ grows, it is enough to understand the 
L-function $L(\rho_t, s)$ for all classical points $t$ near the boundary of $D_p(1)$.  More precisely, we should understand the 
following universal L-function.

\begin{Def} 
The $T$-adic L-function of the tower $K_{\infty}/K$ is the L-function of the $T$-adic character $\rho_T$: 
$$L_{\rho}(T,s):= L(\rho_T, s) = \prod_{x\in U} \frac{1}{1-{\rho_T({\rm{Frob}}_x)}s^{{\rm{deg}}(x)}}\in 1 +s\ZZ_p[[T]][[s]].$$
\end{Def}
If $x$ is a prime of $K$, then the universal character $\rho_T$ is unramified at $x$ if and only if $\rho_T(I_x)=1$. This is true if and only if $I_x =1$. Thus, the ramification 
locus of $\rho_T$ is exactly $P-U$.  
So, the above definition agrees with the usual definition of L-functions. 
The set $U$ is the set of closed points of a smooth affine curve over $\FF_q$. 
Using Crew's work \cite{crew} on the Katz conjecture in the abelian case, 
%or more generally Emerton-Kisin's work \cite{EK} in the general case, 
we deduce the following decomposition of the 
universal L-function, which gives the analytic continuation on the closed unit disk. 

\begin{Thm}\label{Crew} There is a polynomial $P(T, s) \in 1+ \ZZ_p[[T]][s]$ in $s$ and a power series 
$g(T, s) \in 1+(p,T)\ZZ_p[[T]][[s]]$ which converges for all $s\in \ZZ_p[[T]]$  such that 
$$L_{\rho}(T,s) = P(T, s)g(T, s).$$ 
\end{Thm}

As a consequence, we deduce that $g(T,1) \in 1+(p,T)\ZZ_p[[T]]$, 
$$P(T,1) \in \ZZ_p[[T]], \ L(T):= L_{\rho}(T,1)=P(T,1)g(T, 1)\in \ZZ_p[[T]]$$
are well defined power series in $T$ over $\ZZ_p$. Furthermore, for $t\in D_p(1)$ and 
$s\in \CC_p$ with $|s|_p\leq 1$, we have 
$$v_p(L_{\rho}(t, s)) = v_p(P(t,s)).$$ 
 
Recall that for each interior classical point $t=(t_1, \cdots, t_d) \in D^*$, 
the ramification locus of $\rho_t$ is exactly $P-U$.  
Using the definition of L-functions, one checks 
that for $t\in D^*$, we have 
$$L(\rho_{t}, s) = L_{\rho}(T,s)|_{T=t} =L_{\rho}(t, s) $$
In particular, for every interior finite character $\chi=\chi_1\otimes \cdots \otimes \chi_d 
\in X_n^*$, we have $t=(\chi_1(1)-1, \cdots, \chi_d(1)-1) \in D^*$  and thus 
$$L(\chi, s) =L(\rho_t, s) =L_{\rho}(t,s), \ L(\chi, 1) = L_{\rho}(t,1) =L(t).$$
We conclude that 
$$\sum_{\chi \in X_n^*} v_p(L(\chi, 1)) = \sum_{\chi\in X_n^*} v_p(L(\chi_1(1) -1, \cdots, \chi_d(1)-1)).$$ 
Since the set $X^*$ is open (and hence semi-algebraic) in $W^d$, 
Theorem 5.6 in \cite{Mo1} implies that the right side is a polynomial in $p^n$ and 
$n$ of total degree at most $d$, and of degree at most $1$ in $n$. 
Theorem \ref{Mon} is proved. Thus, Greenberg's conjecture is true 
in the function field case when there is no constant subextension. 

We now explain how the above proof implies the result on the $p$-rank of $K_n$. 
Since the character $\rho_T$ is trivial modulo $T$, the L-function $L_{\rho}(T,s)$ modulo $T$ is the same as the zeta function $Z(U,s)$ of $U$. 
This gives the congruence 
$$L_{\rho}(T, s) \equiv Z(K_0, s) \prod_{x\in P-U} (1- s^{{\rm deg}(x)})    = \frac{P(K_0, s)}{(1-s)(1-qs)} \prod_{x\in P-U} (1- s^{{\rm deg}(x)}) \mod T .$$
Replacing $T$ by an interior classical point $t=(t_1, \cdots, t_n)$ in $D^*$, we deduce 
$$L(\rho_t, s)= L_{\rho}(t, s) \equiv \frac{P(K_0, s)}{(1-s)} \prod_{x\in P-U} (1- s^{{\rm deg}(x)}) \mod (t_1, \cdots, t_d).$$ 
Let $\ell_p(t)$ denote the number of $p$-adic unit roots of $L(\rho_t,s)$. 
Comparing the number of $p$-adic unit roots on both sides, one finds that for every interior classical point $t\in D^*$, 
$$\ell_p(t) = r_p(0) -1 +\sum_{x\in P-U} {{\rm deg}(x)}$$
is a constant $c$ indepedent of $t$, where $r_p(0)$ is the $p$-rank of $K_0$ which is the number of $p$-adic unit roots of $P(K_0, s)$.  
Theorem \ref{rank1} is proved. Since the $p$-rank is a geometric invariant, that is, independent of the constant 
extension, it follows that Theorem \ref{rank} is also proved.

\section{The case with constant $\ZZ_p$-extension} 

Assume now that the $\ZZ_p^d$-tower  $K_{\infty}/K$ has a $\ZZ_p$-subtower of constant extension. 
Then, we can write 
$$K_{\infty} =\FF_{q^{p^{\infty}}} L_{\infty}, \ K_n = \FF_{q^{p^n}} L_n,$$
where $L_{\infty}/K$ is a $\ZZ_p^{d-1}$-tower with no constant extension. To avoid triviality, we can assume that $d\geq 2$. 
Replacing $K$ by $K_m$ for some $m$, we may and will assume that the tower $L_{\infty}/K$ satisfies the 
normalization in the previous section and is ramified on the non-empty finite set $P-U$. 
The zeta function of $K_n$ is then 
given by $p^n$-th Adams operation of the zeta function of $L_n$. That is, 
$$Z(K_n, s) = \Phi^{p^n} (Z(L_n, s)) = \frac{\Phi^{p^n}(P(L_n, s))}{(1-s)(1-q^{p^n}s)},$$
where $\Phi^{p^n}$ denotes the $p^n$-th Adams operation. Recall that for a polynomial 
$$H(s) = \prod_{i=1}^m (1-\alpha_i s) = \det(I -A s)\in \ZZ_p[s],$$ 
we have 
$$\Phi^{p^n}(H(s)) = \prod_{i=1}^m (1-\alpha_i^{p^n} s) = \det(I -A^{p^n}s).$$
That is, it raises each reciprocal root (or eigenvalue) to its $p^n$-th power.  
It follows that 
$$\Phi^{p^n}(H(s))|_{s=1} = \det(I -A^{p^n}) =\prod_{\eta^{p^n}=1} \det(I-\eta A) =\prod_{\eta^{p^n}=1} H(\eta).$$
Apply the previous interior decomposition proof in the no constant extension case to the $\ZZ_p^{d-1}$-tower $L_{\infty}/K$, we deduce that 
for each subset $S\subset P-U$,  
there is a polynomials $P_S(T_S,s)\in 1 +s\ZZ_p[[T_S]][s]$ such that 
$$v_p(\frac{h_n}{h_0}) = \sum_{S\not= \phi, S \subset P-U} \Big( \sum_{\chi \in X_{n,S}^*}\sum_{\eta^{p^n}=1}
 v_p(P_S(\chi_1-1, \cdots, \chi_{d-1}(1)-1, (\eta -1)+1))\Big).$$
The set $X_S^* \times W$ is an open subset of $W^{d(S)}\times W = W^{d(S)+1}$. By Theorem 5.6 in \cite{Mo1} again, 
the right side is a polynomial in $p^n$ and $n$ of total degree at most $d(S)+1\leq d$, and of degree at most $1$ in $n$. 
The theorem is proved.  

\section{Stability of zeta functions in $p$-adic Lie towers}

In this final section,  $K$ is a global function field of characteristic $p$ with constant field $\FF_q$. 
We propose several conjectures on possible stability of zeta functions in various $p$-adic Lie towers over $K$, 
vastly extending the function field Greenberg conjecture  in several different directions. 

Let $G_K ={\rm Gal}(K^{sep}/K)$ denote 
the absolute Galois group of $K$. Let 
$$\rho: G_K \longrightarrow {\rm GL}_d(\ZZ_p)$$ 
be a continuous $p$-adic representation of $G_K$ with infinite image, ramified at finite number of primes of $K$. The image $G_{\infty}$ of this representation is a compact $p$-adic Lie group of dimension ${\rm dim}(G)>0$.  
The fixed field $K_{\infty}$ of the kernel $Ker(\rho)$ is a Galois extension of $K$ with Galois group isomorphic to $G_{\infty}$.  
That is, ${\rm Gal}(K_{\infty}/K) \cong G_{\infty}$. For integer $n\geq 0$, the kernel of the reduction 
$$\rho_n: G_K  \longrightarrow {\rm GL}_d(\ZZ_p) \longrightarrow  {\rm GL}_d(\ZZ_p/p^n\ZZ_p)$$ 
produces a Galois extension $K_n$ of $K$ whose Galois group is isomorphic to the image $G_n$ of $\rho_n$. 
This produces a $p$-adic Lie tower
$$K_{\infty}=\bigcup_{n=0}^{\infty} K_n \supset \cdots  \supset K_n \supset \cdots \supset K_1 \supset K_0=K, 
\ {\rm Gal}(K_{n}/K) \cong G_{n},$$
which ramifies at finite number of primes of $K$. 
In the spirit of Iwasawa theory, our basic problem is to understand 
possible stable behavior of various arithmetic properties of this tower of global fields $K_n$ as $n$ grows. 
This naturally raises many interesting open problems. 
We state a few of them as conjectures below. 
The first one concerns the $p$-ranks and class numbers of this tower of global fields.

\begin{Con}\label{Conj1}  Let $K_{\infty}/K$ be a $p$-adic Lie extension as above. 
Let $r_p(n)$ denote the $p$-rank of $K_n$ and let $h_n$ denote the class number of $K_n$. 

\item{(1)} ($p$-Rank stability). There is a polynomial $R(x)\in  \QQ[x]$ of degree at most ${\rm dim}(G)$ depending on the tower such that for all sufficiently large $n$, 
we have 
$$r_p(n) = R(p^n).$$

\item{(2)}($p$-Class number stability). There is a polynomial $E(x,y)\in  \QQ[x,y]$ of total degree at most ${\rm dim}(G)$ depending on the tower such that for all sufficiently large $n$, 
we have 
$$v_p(h_n) = E(p^n, n).$$
\end{Con}

Our previous results for $\ZZ_p^d$-towers imply that this conjecture is true when the $p$-adic Lie group $G_{\infty}$ is abelian. 
The conjecture can be viewed as an attempt to extend geometric Iwaswa theory to non-abelian extensions. 
Part (1) of the conjecture can be proved if each ramified prime is totally ramified.  This provides a positive 
evidence to part (2) of the conjecture which seems significantly more difficult. 

For a given $K$, there are too many  $p$-adic Lie towers, most of them are not natural. 
To state further stability conjectures, we need to assume that the tower $K_{\infty}/K$ is natural in some sense, 
that is, arising from algebraic geometry. The tower $K_{\infty}/K$ 
is called algebraic geometric if the corresponding Galois representation $\rho$ arises 
from a $p$-adic \'etale cohomology of a smooth proper variety $X$ defined over the global field $K$. 
The algebraic geometric tower $K_{\infty}/K$ is further called ordinary if the variety $X$ is generically ordinary over $K$. 
An algebraic geometric representation automatically ramifies at finite number of primes of $K$. 
In the case that $X=A$ is an ordinary abelian variety of dimension $g$ over $K$,  the first $p$-adic \'etale cohomology 
can be explicitly constructed from the $p$-adic Tate module: 
$$T_p(A) = \lim_{\leftarrow~k} A [p^k],$$
which gives a $g$-dimensional $p$-adic representation of $G_K$. This example 
is extremely important, already in the case $g=1$, because of its close connection to 
$p$-adic automorphic forms.

It was conjectured by Dwork \cite{Dw} and proved by the author \cite{Wa1}\cite{Wa2} that the L-function $L(\rho, s)$ of any algebraic geometric $p$-adic 
representation $\rho$ is $p$-adic meromorphic everywhere in $s\in \CC_p$. This nice analytic property suggests that the $p$-adic Lie tower 
arising from an algebraic geometric $p$-adic representation should have good stable geometric and arithmetic properties. 
The simplest is the following genus stability conjecture. 

\begin{Con}[Genus stability]\label{genus} Let $g_n$ denote the genus of $K_n$.  Assume that the tower $K_{\infty}/K$ 
is algebraic geometric and ordinary. Then there is a polynomial $G(x) \in \QQ[x]$ of degree at most ${\rm dim}(G)+1$ depending on the tower such that for all sufficiently large $n$, 
we have 
$$ g_n=G(p^n).$$
\end{Con}

A major progress toward this conjecture has been made by Joe Kramer-Miller in \cite{Kr} and in his forthcoming work. 
Motivated by this conjecture, the number field analogue on discriminant stability has been proved recently by James Upton \cite{Up}.   
If the algebraic geometric tower is non-ordinary, then a slightly weaker genus stability is still expected to hold as suggested 
by Krammer-Miller. In this more general case, one would need several polynomials $G_j(x) \in \QQ[x]$ ($1\leq j \leq m$) to express 
the genus $g_n$ such that $g_{mn+j} = G_j(n)$ for all large $n$.

Next, we move to the deeper slope stability for algebraic geometric $p$-adic Lie towers. Without loss of essential generality, we can 
and will assume that the tower has no proper constant subextension below.  

The zeta function 
$Z(K_n, s)$ of $K_n$  is a rational function in $s$ of the form 
$$Z(K_n, s) = \frac{P(K_n, s)}{(1-s)(1-qs)}, \ P(K_n, s) \in 1 +s\ZZ[s],$$
where the zeta polynomial $P(K_n, s)$ is a polynomial of degree $2g_n$ and $g_n= g(K_n)$ denotes the 
genus of $K_n$. Recall that the class number is given by the special zeta value $h_n = P(K_n, 1)$ 
and the the $p$-rank $r_p(n)$ is the number of $p$-adic unit roots of $P(K_n, s)$. 
Thus, the above two conjectures are really about the $p$-adic stability of some aspects of the zeta function 
$Z(K_n, s)$ as $n$ grows. 

A deeper property is the possible stability of higher slopes for the zeta function, not just the 
$p$-rank which is the slope zero part. To describe this, 
we write 
$$P(K_n, s) = \prod_{i=1}^{2g_n} (1-\alpha_i(n)s) \in \CC_p [s], \ \  0\leq v_q(\alpha_1(n)) \leq \cdots \leq v_q(\alpha_{2g_n}(n)) \leq 1, $$
where the $q$-adic valuation is normalized such that $v_q(q)=1$. These rational numbers $\{v_q(\alpha_1(n)), \cdots, v_q(\alpha_{2g_n}(n)) \}$ 
are called the $q$-slopes of $K_n$. They are symmetric in the interval $[0,1]$ by the functional equation.  
For fixed rational number $\alpha \in [0,1]$, let $r_p(n, \alpha)$ denote the multiplicity of the 
slope $\alpha$ in the slope sequence $\{v_q(\alpha_1(n)), \cdots, v_q(\alpha_{2g_n}(n)) \}$. For $\alpha=0$, $r_p(n, 0) =r_p(n)$ is 
simply the $p$-rank of $K_n$. The integer $r_p(n, \alpha)$ is called the slope $\alpha$-rank of $K_n$.  
We would like to understand how the slope $\alpha$-rank $r_p(n, \alpha)$ and more generally how the full slope sequence vary when $n$ grows. 
The following conjecture has three parts with increasing level of difficulty. 

\begin{Con}\label{slope}
Assume that the $p$-adic Lie tower $K_{\infty}/K$ is algebraic geometric and ordinary with no proper constant subextension. 

\item{(1)} (Slope uniformity). The $q$-slopes  
$$\{ v_q(\alpha_1(n)), \cdots, v_q(\alpha_{2g_n}(n))\} \subset [0, 1] \cap \QQ \subset [0, 1]$$ 
are equi-distributed in the interval $[0, 1]$ as $n$ goes to infinity. 

\item{(2)} (Slope $\alpha$-rank stability). For each fixed rational number $\alpha \in [0, 1]$, 
there is a  polynomial $R_{\alpha}(x) \in \QQ[x]$ of degree at most ${\rm dim}(G)$ depending on the tower such that for all sufficiently large $n$, 
we have 
$$ r_p(n, \alpha)=R_{\alpha}(p^n).$$

\item{(3)} (Slope stability). There is a positive integer $n_0$ depending on the tower such that the re-scaled $q$-slopes
$\{ p^nv_q(\alpha_1(n)), \cdots, p^nv_q(\alpha_{2g_n}(n))\}$ for all $n> n_0$ are determined explicitly by 
their values for $0\leq n \leq n_0$, using a finite number of arithmetic progressions. 
\end{Con}

If the tower is non-ordinary, one expects similar and possibly slightly weaker stability results.  If the 
tower does not come from algebraic geometry, the conjecture is false in general.  
A remarkable recent work of Kosters-Zhu \cite{KZ} suggests that the genus stability conjecture always implies the slope uniformity 
conjecture. In fact, they have proved this implication when $K=\FF_q(x)$ is the rational function field and $G_{\infty}=\ZZ_p$ in \cite{KZ} and in their forthcoming work.  
We do not have a precise formulation about the harder slope stability conjecture in such a generality yet. 
The essence is a finiteness property that all the slopes  for large $n$ are determined by a finite number of arithmetic progressions. 
Needless to say that this conjecture (and even the weaker $\alpha$-rank stability conjecture) is wide open, 
but interesting progress has been made in various special cases when the 
$p$-adic tower $K_{\infty}/K$ is abelian, see \cite{DWX}\cite{RWXY}\cite{Li}\cite{KZ}. In the case that the tower is the Igusa $\ZZ_p^*$-tower over a modular curve, the 
above slope stability conjecture is partially related to the ghost conjecture \cite{BP} and the spectral halo conjecture on $p$-adic modular forms, see \cite{liu-wan-xiao}\cite{wan-xiao-zhang} for recent progress and further references.

%\vfill
%\eject


\begin{thebibliography}{9999}

%\bibitem[BG]{buzzard-gee} K. Buzzard and T. Gee, Slopes of modular forms,  	{\tt arXiv:1502.02518}. 


\bibitem[BP]{BP}J. Bergdall and R. Pollack,  Slopes of modular forms and the ghost conjecture, 
IMRN.,  to appear, {\tt arXiv:1607.04658}. 

%\bibitem[CM]{coleman-mazur}
%R. Coleman and B. Mazur,
%The eigencurve, in {\it Galois representations in arithmetic algebraic geometry}, 
%{Cambridge Univ. Press, Cambridge}, 1998, 1-113. 


\bibitem[Cr]{crew} R. Crew,  $L$-functions of $p$-adic characters and geometric Iwasawa theory, 
{\it Invent. Math.},  {\bf 88} (1987), no. 2, 395-403.


\bibitem[Cu] {Cu}A. Cuoco,  The growth of Iwasawa invariants in a family, 
{\it Comp. Math.},  41(1980), pp. 415-417. 


\bibitem[CM]{CM} A. Cuoco  and P. Monsky, Class numbers in $\ZZ_p^d$-extensions, 
{\it Math. Ann.},  255(1981), pp. 235-258. 

\bibitem[DWX]{DWX} C. Davis, D. Wan and L. Xiao, Newton slopes for Artin-Schreier-Witt towers, 
{\it Math. Ann.}, {\bf 364} (2016), no. 3, 1451-1468.   


\bibitem[Dw]{Dw} B. Dwork, Normalized period matrices II, {\it Ann. Math.}, {\bf 98} (1973), 1-57. 


%\bibitem[EK]{EK} M. Emerton and M. Kisin, Unit L-functions and a conjecture of Katz, 
%{\it Ann. of Math.}.  (2) 153 (2001), no. 2, 329-354. 

%\bibitem[Gk]{EG} E. Grosse-Kl\"onne, On families of pure slope L-functions, 
%{\it Documenta Math.}, 8 (2003) 1--42. 

\bibitem[GK]{GK} R. Gold and H. Kisilevsky, On geometric $\ZZ_p$-extensions of function fields, 
{\it Manuscripta Math.}, {\bf 62} (1988), 145--161. 

\bibitem[Gr]{Gr} R. Greenberg, The Iwasawa  invariants of $\Gamma$-extensions of a fixed number field, 
{\it Amer. J.  Math.}, Vol. 95, No. 1 (1973), 204-214. 

%\bibitem[Ha]{Ha} D. Haessig,  $L$-functions of symmetric powers of Kloosterman sums (unit root L-functions and $p$-adic estimates), 
%Math. Ann. (2016). doi:10.1007/s00208-016-1454-6. 

\bibitem[Ka]{Ka} N. Katz,  Travaux de Dwork, S\'eminaire Bourbaki, expos\'e 409(1971/72), 
Lecture Notes in Math., 317(1973), 167-200. 


\bibitem[KW]{KW}M. Kosters and D. Wan,  Genus growth in $\ZZ_p$-towers of function fields, 
{\it Proc. Amer. Math. Soc.},  146 (2018), no. 4, 1481-1494. 

\bibitem[KZ]{KZ}M. Kosters and J. Zhu, Slopes of L-functions in genus stable $\ZZ_p$-covers of the projective line, 
{\it J. Number Theory},  187 (2018), 430-452. 

\bibitem[Kr]{Kr}J. Kramer-Miller, The monodromy of $F$-isocrystals with logarithmic decay, arXiv:1612.01164. 

\bibitem[Li]{Li}X. Li,  The stable property of Newton slopes for polynomial Witt towers,  
{\it J. Number Theory},  185 (2018), 144-159. 

\bibitem[LZ]{LZ}C. Li and J. Zhao, Iwasawa theory of $\ZZ^d_p$ -extensions over
global function fields,  {\it Exposition. Math.},  15 (1997), no. 4, 315-337. 



\bibitem[LW]{LW} C.  Liu and D. Wan,
$T$-adic exponential sums over finite fields, {\it Algebra Number Theory},  {\bf3} (2009), no. 5, 489--509. 


%\bibitem[LLN]{LLN} C. Liu, W. Liu,  and C.  Niu,  Generic T-adic exponential sums in one variable, 
%{\it J. Number Theory}, {\bf 166} (2016), 276-297. 


%\bibitem[LWei]{liu-wei} C. Liu and D. Wei, 
%The L-functions of Witt coverings, {\it Math. Z.}  {\bf 255} (2007), 95--115. 

\bibitem[LWX]{liu-wan-xiao}
R. Liu, D. Wan, and L. Xiao,
Slopes of eigencurves over the boundary of the weight space, 
{\it Duke Math. J.}, 166(2017), no. 9, 1739-1787. 

\bibitem[MW]{MW} B. Mazur and A. Wiles, Analogies between function fields and number fields, 
{\it Amer. J. Math.},  105 (1983), no. 2, 507-521. 


\bibitem[Mo1]{Mo1} P. Monsky, On $p$-adic power series, 
{\it  Math. Ann.},  255(1981), pp. 217-227. 

\bibitem[Mo2]{Mo2} P. Monsky, Som invariants of $\ZZ_p^d$-extensions, 
{\it  Math. Ann.},  255(1981), pp. 229-233. 

\bibitem[Mo3]{Mo3} P. Monsky, $p$-Ranks of class groups in $\ZZ_p^d$-extensions, 
{\it  Math. Ann.},  263(1983), pp. 509-514. 

\bibitem[Mo4]{Mo4} P. Monsky, Fine estimates for the growth of $e_n$ in $\ZZ^d_p$-extensions, 
in: Algebraic Number Theory| in honor of K. Iwasawa, J. Coates, R. Greenberg, B.
Mazur and I. Satake (eds.), Advanced Studies in Pure Math. 17(1989), pp. 309-330. 

%\bibitem[Mo]{Mo} P. Monsky, Formal cohomology III, {\it Ann. Math.}, {\bf 93} (1971), 315-343. 


%\bibitem[OY]{OY} Y. Ouyang  and J. Yang,  Newton polygons of L functions of polynomials $x^d + ax$, 
%{\it J. Number Theory}, {\bf 160} (2016), 478-491. 

%\bibitem[W2]{wan4}
%D. Wan, Dimension variation of classical and $p$-adic modular forms, 
%{\it Invent. Math.} {\bf 133} (1998), no. 2, 449--463. 

\bibitem[RWXY]{RWXY} R. Ren, D. Wan, L. Xiao and M. Yu, Slopes for higher rank Artin-Schreier-Witt towers, 
{\it Trans. Amer. Math. Soc.},  to appear, {\tt arXiv:1605.02254}. 

%\bibitem[Tate]{Tate} {J.T. Tate}, 
%$p$-divisible groups, {\it Proc.  Conf.  Local  Fields} (Driebergen, 1966), Springer, Berlin, 1967.

%\bibitem[W96]{wan96} D. Wan, Meromorphic continuation of $L$-functions of $p$-adic representations, {\it Ann. Math.}, {\bf 143}(1996), 469-498. 

\bibitem[Up]{Up} J. Upton,  Discriminant-stability in $p$-adic analytic towers of number fields, arxiv.org/abs/1801.03056. 


\bibitem[Wa0]{Wa0} D. Wan, Meromorphic continuation of $L$-functions of $p$-adic representations, 
{\it Ann. Math.}, 143( 1996), 469-498. 

\bibitem[Wa1]{Wa1} D. Wan, Rank one case of Dwork's conjecture, 
{\it J. Amer. Math. Soc}. {\bf 13} (2000), no. 4, 853-908. 

\bibitem[Wa2]{Wa2} D. Wan, Higher rank one case of Dwork's conjecture, 
{\it J. Amer. Math. Soc}. {\bf 13} (2000), no. 4, 807-852. 

\bibitem[WXZ]{wan-xiao-zhang} D. Wan, L, Xiao and J.  Zhang,  Slopes of eigencurves over boundary disks, 
{\it  Math. Ann.}, 369(2017), no. 1-2, 487-537. 


%\bibitem[W99]{wan2}
%D. Wan, Dwork's conjecture on unit root zeta functions, 
%{\it Ann.  Math.} {\bf 150} (1999), no. 3, 867--927. 

%\bibitem[W4]{wan3}
%D. Wan, Variation of $p$-adic Newton polygons for L-functions of exponential sums, 
%{\it Asian J. Math.} {\bf  8} (2004), no. 3, 427--471.

 


\end{thebibliography}
\end{document}